\documentclass[11pt,reqno]{amsart}

\newcommand{\be}{\begin{eqnarray}}
\newcommand{\ee}{\end{eqnarray}}
\def\es{\emptyset}

\def\ld{\stackrel{{\rm LD}}{\longrightarrow}}

\def\Sf{{\sf S}}
\def\Lb{{\bf \Lambda}}
\def\Gamb{{\bf \Gamma}}
\def\Gam{\Gamma}
\def\Xb{{\bf X}}
\def\be{{\bf e}}
\newcommand{\supp}{\mbox{\rm supp}}

\newcommand{\R}{{\mathbb R}}

\newcommand{\Z}{{\mathbb Z}}
\newcommand{\Nat}{{\mathbb N}}

\newcommand{\Fk}{{\mathcal F}}

\newcommand{\Lk}{{\mathcal L}}

\newcommand{\Sk}{{\mathcal S}}
\newcommand{\Tk}{{\mathcal T}}

\newcommand{\dist}{\mbox{\rm dist}}

\newcommand{\lam}{\lambda}
\def\om{\omega}
\def\gam{\gamma}
\def\Lam{\Lambda}

\newcommand{\Int}{{\rm int}}

\newcommand{\diam}{{\rm diam}}

\newtheorem{theorem}{Theorem}[section]
\newtheorem{lemma}[theorem]{Lemma}

\newtheorem{prop}[theorem]{Proposition}

\theoremstyle{definition}
\newtheorem{defi}[theorem]{Definition}

\theoremstyle{remark}

\numberwithin{equation}{section}

\begin{document}

\title[Pseudo-self-affine tilings in $\R^d$]{Pseudo-self-affine tilings in $\R^d$}
\date{\today}
\author{Boris Solomyak} 
\thanks{
supported in part by NSF grant \#DMS-0355187.}

\address{Boris Solomyak, Box 354350, Department of Mathematics,
University of Washington, Seattle, WA 98195, USA,}
\email{solomyak@math.washington.edu}

\begin{abstract}
It is proved that every
pseudo-self-affine tiling in $\R^d$ is mutually locally derivable with
a self-affine tiling. A characterization of pseudo-self-similar tilings in
terms of derived Vorono\"{\i} tessellations is a corollary.
Previously, these results were obtained in the planar case, jointly
with Priebe Frank. The new approach is based on the theory of
graph-directed iterated function systems and substitution Delone sets
developed by Lagarias and Wang.
\end{abstract}

\maketitle

\section{Introduction}

\thispagestyle{empty}

Self-affine tilings have been much studied, see \cite{SolTil,Robi} and
the references therein.
They arise, in particular, in connection with Markov partitions of
toral automorphisms and as models for quasicrystals.
A {\em self-affine tiling} has the property that, if we expand it by
a certain expanding linear map, 
then the original tiling may be obtained by subdividing
the expanded tiles according to a prescribed rule. 

A {\em pseudo-self-affine tiling} has 
a less rigid hierarchical structure: 
if the entire tiling is expanded by a certain
linear map, then one can recover the original tiling at
any point by looking in a finite ``window'' around that point in the
expanded tiling. If the expanding linear map is a similitude, we say that
the tiling is {\em pseudo-self-similar}.
A famous example of a pseudo-self-similar tiling is the
Penrose tiling with rhombi tiles.

Pseudo-self-affine tilings were introduced by 
N. Priebe Frank \cite{Thesis,PF}.
E. A. Robinson, Jr. conjectured that such tilings are {\em mutually locally
derivable} from self-affine tilings (precise definitions are given in the
next section).
In \cite{PrSo}, joint with N. Priebe Frank, 
we settled the conjecture for pseudo-self-similar tilings
in $\R^2$ and used it to complete the characterization of
pseudo-self-similar tilings
in $\R^2$ started in \cite{PF}.

Here we establish the conjecture in full generality.
The only caveat is that
the tiles of the resulting self-affine tiling need not be connected
(in \cite{PrSo} the tiles were topological disks);
however, having such a tiling is sufficient for many purposes; in particular,
it yields a characterization of pseudo-self-similar tilings in $\R^d$.
The approach is different: whereas in \cite{PrSo} we used
the method of ``redrawing the boundary,'' here we obtain
the prototiles as an attractor of a graph-directed iterated function system,
applying the results of Lagarias and Wang \cite{LaWa}.

The work of A. M. Vershik and co-authors on arithmetic constructions
of Markov and sofic partitions \cite{Ver,KV,SiVer} 
has many points of contact with the theory of self-affine tilings; it played
an important role in the development of the subject.
I am happy to dedicate this paper to Anatoly Moiseevich Vershik,
with gratitude and affection.


\section{Preliminaries}

We begin with tiling preliminaries, following \cite{PrSo} and \cite{LMS}.
See also \cite{Robi} for a recent survey.

We fix a set of types (or colors) labeled by $\{1,\ldots,m\}$. 
A {\em tile} in $\R^d$ is defined as a pair $T = (A,i)$ where 
$A  = \supp(T)$ (the support of $T$) is a compact set in $\R^d$ which is the
closure of its interior, and $i = \ell(T)\in \{1,\ldots,m\}$ is the type of $T$.
We emphasize that the tiles are not assumed to
be homeomorphic to the ball or even connected.
A {\em tiling} of $\R^d$ is a set $\Tk$ of tiles such that $\R^d = \bigcup
\{\supp(T):\ T\in \Tk\}$ and distinct tiles have disjoint interiors.

A $\Tk$-{\em patch} $P$ is a finite subset of the tiling $\Tk$.
Denote by $\Tk^*$ the set of all $\Tk$-patches.
The {\em support of a patch} $P$ is defined by
$\supp (P) = \bigcup\{\supp(T):\ T\in P\}$. The {\em diameter of a patch}
$P$ is $\diam(P)=\diam(\supp(P))$. The {\em translate} of a tile $T=(A,i)$
by a vector $g\in \R^d$ is $T+g = (A+g, i)$. The translate
of a patch $P$ is $P+g = \{T+g:\ T\in P\}$. We say that two patches
$P_1,P_2$ are {\em translationally equivalent} if $P_2 = P_1+g$ for some
$g\in \R^d$. 

A {\em patch with a marked tile} is a pair $(P,T)$ where $P$ is a patch and
$T\in P$. Two patches with marked tiles $(P_1,T_1)$ and $(P_2,T_2)$ are
said to be translationally equivalent if $P_2=P_1+g$ and $T_2 = T_1+g$ for
some $g\in \R^d$.

Given an invertible linear map $\psi$ in $\R^d$ and a tiling $\Tk$ 
we can consider a new tiling
$\psi\Tk = \{(\psi(\supp(T)),\ell(T)):\ T\in \Tk\}$.

For a tile $T$ and $\Lam \subset \R^d$ we denote
$$
T + \Lam:= \{T + g:\ g\in \Lam\}.
$$

We always assume that:
\begin{itemize}
\item Any two $\Tk$-tiles with the same type (color)
are translationally equivalent.
(Hence there are finitely many $\Tk$-tiles up to translation.)
\item the tiling $\Tk$ has {\em finite local complexity} (FLC), that is,
for any $R>0$ there are finitely many $\Tk$-patches of diameter less
than $R$ up to translation equivalence.
\end{itemize}

Given a tiling $\Tk$, we can choose one tile for each translation
equivalence class, thereby obtaining a 
{\em prototile set} $\{T_1,\ldots,T_m\}$ (note that $T_i \in \Tk$ by
assumption).
Then we can write
\begin{equation} \label{eq-delon1}
\Tk = \bigcup_{i=1}^m (T_i + \Lam_i)
\end{equation}
for some uniformly discrete sets $\Lam_i$.

\medskip

\noindent {\bf Notation.} We fix a metric in $\R^d$ (equivalent to the
Euclidean metric), and write, for $F \subset \R^d$:
$$
N_R(F) = \{x\in \R^d:\ \dist(x,F) \le R\}.
$$
Thus $N_R(x) = N_R(\{x\})$ is the closed ball of radius $R$ centered at $x$.

\begin{defi} \label{def-rep}
A tiling $\Tk$ is called {\em repetitive} if for any patch $P\subset \Tk$
there is a real number $R>0$ such that for any $x\in \R^d$ there is a
$\Tk$-patch $P'$ such that $\supp(P')\subset N_R(x)$ and $P'$ is a
translate of $P$.
The minimal such $R$, denoted $R(P)$, is called the {\em repetitivity radius}
of $P$.
\end{defi}

If the tiling $\Tk$ is repetitive, then the sets $\Lam_i$ in (\ref{eq-delon1})
are both uniformly discrete and relatively dense. Such sets are called
{\em Delone sets}.

\medskip

\noindent {\bf Tiling Dynamical Systems.} Although in this paper 
we do not deal with
dynamical systems directly, they provide a useful framework
for the concepts and results. We recall the set-up briefly.

The {\em tiling space} is 
$X_{\Tk} = \overline{\{-g+\Tk : g \in \R^d\}}$, where $X_{\Tk}$ carries a
well-known tiling topology. It is based on the idea that two tilings are
close if after a small translation they agree on a large ball around
the origin. The reader is referred to \cite{Robi} for details.
The space $X_\Tk$ is compact \cite{Rud} (see also \cite{Robi}), 
and the group $\R^d$ acts
on it continuously by translations, so that we get a
{\em tiling dynamical system} associated with $\Tk$.
It is well-known (see \cite{Robi}) that $\Tk$ is repetitive if and only if
the tiling dynamical system $(X_{\Tk},\R^d)$ is minimal, that is, its every
orbit is dense.

\medskip

Next we define local derivability, which is the key concept in this paper.
First, some more notation: to
a subset $F\subset \R^d$ and a tiling $\Tk$ we associate a
$\Tk$-patch as follows:
$$
[F]^{\Tk} = \{T\in \Tk:\ \supp(T)\cap F \ne \emptyset\}.
$$
We write $[x]^{\Tk} = [\{x\}]^{\Tk}$ for notational convenience.

\begin{defi} (See \cite{Baake-ld}.) \label{def-ld}
Let $\Tk_1$ and $\Tk_2$ be two tilings. We say that $\Tk_2$ is
{\em locally derivable (LD) from} $\Tk_1$ with a radius $R>0$
if for all $x,y \in \R^d$,
\begin{equation} \label{ld} 
[N_R(x)]^{\Tk_1} = [N_R(y)]^{\Tk_1} + (x-y) \ \Rightarrow \
[x]^{\Tk_2} = [y]^{\Tk_2} + (x-y).
\end{equation}
We write $\Tk_1 \ld \Tk_2$ to denote that $\Tk_2$ is LD from $\Tk_1$.
If $\Tk_1\ld\Tk_2$ and $\Tk_2\ld\Tk_1$, then we say that
$\Tk_1$ and $\Tk_2$ are {\em mutually locally derivable (MLD)}.
\end{defi}

\noindent {\bf Remark.}
It is clear that MLD is an equivalence relation.
If $\Tk_1 \ld \Tk_2$ then there is a factor map of the corresponding
topological dynamical systems, see \cite{PF}. It is the tiling analog
of a finite block code in symbolic dynamics. If $\Tk_1$ and $\Tk_2$ are
MLD, then the associated tiling dynamical systems are topologically
conjugate, but the converse is not true \cite{petersen,rasad}.

\medskip

Let $\phi:\ \R^d \rightarrow \R^d$
be an expanding linear map, that is, all its eigenvalues are greater than
1 in modulus. Then there is a norm $|\cdot|$ and $\lambda>1$ such that
\begin{equation} \label{expans}
|\phi g|\ge \lambda |g|,\ \ \ g\in \R^d.
\end{equation}
We fix such a norm, and use the corresponding metric, whenever we have
an expanding linear map.
(Alternatively, we can pass from $\phi$ to $\phi^\ell$ for appropriate
$\ell\in \Nat$, and use the Euclidean norm.)

\begin{defi} \label{def-hier}
Let $\phi:\ \R^d\to \R^d$ be an expanding linear map.
A repetitive FLC tiling $\Tk$
is called a {\em pseudo-self-affine tiling with expansion $\phi$} if
$\phi \Tk \ld \Tk$.

A repetitive FLC tiling $\Tk$ is called
a {\em self-affine tiling with expansion $\phi$} if

{\bf (i)} for any tile $T = (A,i)\in \Tk$, there is a $\Tk$-patch 
$\om(T)$ such that $\supp(\om(T)) = \phi A$;

{\bf (ii)} for any tile $T$ and $g \in \R^d$,
$$
(T\in \Tk,\ T+g \in \Tk)\ \Rightarrow\
\om(T+g) = \om(T) + \phi g.
$$
\end{defi}

It is easy to see that a self-affine tiling
is pseudo-self-affine. Also note that the property of being pseudo-self-affine
is preserved under MLD.


\section{New results} 
                                                                                
\begin{theorem} \label{thm-mld}
Let $\Tk$ be a pseudo-self-affine tiling of $\R^d$ with expansion $\phi$.
Then for any $k\in \Nat$ sufficiently large, there exists
a tiling $\Tk'$ which is self-affine with
expansion $\phi^k$, such that $\Tk$ is MLD with $\Tk'$.
\end{theorem}
                                                                                
It is perhaps possible to make sure that the tiles of $\Tk'$ are connected
and even homeomorphic to a ball,
but this would require additional work.

Theorem~\ref{thm-mld} implies that the dynamical systems corresponding
to pseudo-self-affine tilings have the same properties as those
corresponding to self-affine tilings. In particular, they are uniquely
ergodic, not strongly mixing; there are sufficient
conditions for weak mixing, etc.,
see \cite{SolTil}.

\medskip

\noindent {\bf Remark.} 
One can consider tilings that are not translationally
finite (but have FLC with respect to
a larger group, e.g., with respect to all orientation-preserving
Euclidean isometries), such as the ``pinwheel tiling'' and its relatives,
see \cite{HRS} and references there. Recently, B. Rand \cite{rand}
generalized the planar result of \cite{PrSo} to this setting.
It would be interesting to know if this can also be done in higher 
dimensions.

\medskip

A pseudo-self-affine tiling with expansion $\phi$ is said to be
{\em pseudo-self-similar} if $\phi$ is a similitude, i.e., if
$|\phi(x)-\phi(y)| = \lam|x-y|$ for all $x,y\in \R^d$ for some $\lam >1$.
Theorem~\ref{thm-mld}  allows us to complete the characterization
of pseudo-self-affine tilings in $\R^d$ for $d\ge 3$ (the planar case
was done in \cite{PrSo}). 

\begin{defi} (See \cite{PF}.) \label{def-dvor} 
Suppose that $\Tk$ is a repetitive tiling of $\R^d$.
Let $r>0,\  P_r = [B_r(0)]^{\Tk}$ and create the {\em locator set}
$$
\Lk_r = \{q \in \R^d \text{ such that there exists } P \subset \Tk
\text{ with } P_r = P - q\}.
$$
Let $R(r)$ be the repetitivity
radius of $P_r$ so that every ball of radius $R(r)$ in $\Tk$ contains
a translate of $P_r$. The
{\em derived Vorono\"{\i} tiling} $\Tk_r$ has a tile $t_q$ for each
$q \in \Lk_r$ with support
$$
\supp(t_q) = \{x \in \R^d:\ |q - x| \le |q' - x| \text{ for all }
q' \in \Lk_r\};
$$
$t_q$ is labeled by the translational equivalence class of the patch
$[B_{2R(r)}(q)]^{\Tk}$.
The {\em derived Vorono\"{\i} family} is defined by
$\Fk(\Tk) = \{\Tk_r:\ r>0 \}$.
                                                                                
Given an expanding similitude $\phi:\R^d \rightarrow \R^d$, we say a
family of tilings $\Fk$ is {\em $\phi$-finite} if there exist $\Sk_1, \ldots,
\Sk_M$ in $\Fk$ so that for any $\Tk \in \Fk$, there is an $i \in \{1,\ldots,
M\}$ and  $j \in \Z^+$ with $\Tk = \phi^j\Sk_i$.
(Here we identify two tilings if they are equal up to a one-to-one
correspondence between the label sets.)
\end{defi}

\begin{theorem} \label{finite.dv}
A non-periodic, repetitive tiling of $\R^d$ is
pseudo-self-similar if and only if its derived Vorono\"{\i} family is
$\psi$-finite for an expanding, orientation-preserving similitude $\psi$.
\end{theorem}

The reader is referred to \cite[Sect.\,6]{PrSo} where the proof is
given in the 2-dimensional case. In fact, the theorem about a 
pseudo-self-similar tiling being MLD to a self-similar one 
was the only place where dimension 2 was used in \cite{PrSo}.
In view of Theorem~\ref{thm-mld}, the result now
transfers to $\R^d$ for arbitrary $d$.


\section{Proof of Theorem~\ref{thm-mld}}

The main difference with \cite{PrSo} is that there we applied a map to ``redraw
the boundaries'' of tiles, whereas here we apply a map to the 
tiles themselves. The
prototiles of the self-affine tiling will be obtained as
attractors of a graph-directed iterated function system
(a solution of a system of set equations). The proof then proceeds
via the theory of substitution Delone sets developed by Lagarias and
Wang \cite{LaWa}.

\subsection{Reduction.}
We are going to ``recode'' the tiling keeping the supports unchanged, but
increasing the number of labels (this is similar to the
``higher block presentation'' in Symbolic Dynamics, see \cite{LM}).
The new label of a $\Tk$-tile $T$ will be the
equivalence class of the patch $[N_L(\supp(T))]^{\Tk}$ with the marked tile
$T$, for some $L>0$. This allows us to prove the following
                                                                                
\begin{prop} \label{prop-reduce}
Let $\Sk$ be a pseudo-self-affine tiling with expansion $\phi$.
Then there exists a tiling $\Tk$ which is MLD with $\Sk$, has the same
tile supports as $\Sk$ (differs only in labels), such that
for any $T,T' \in \Tk$, $g\in \R^d$, for any $k\ge 1$,
\begin{equation} \label{eq1}
T' = T + g \ \Rightarrow\ [\supp(T')]^{\phi^{-k} \Tk}
= [\supp(T)]^{\phi^{-k} \Tk} + g.
\end{equation}
\end{prop}
                                                                                
The proof of the proposition is straightforward but lengthy,
so we postpone it to the next section.
In view of the proposition, we can assume that the pseudo-self-affine tiling
$\Tk$ in Theorem~\ref{thm-mld} satisfies (\ref{eq1}).

Let
$d_M=d_M(\Tk) = \sup\{\diam(\supp(T)):\ T\in \Tk\}$, which is finite by finite
local complexity. Similarly, we can find $\eta = \eta(\Tk)>0$ such that
every tile support contains a closed ball of radius $\eta$ in its interior.
Recall that $\lam>1$ is the lower bound on the expansion under
the linear map $\phi$.
Fix $k\in \Nat$ such that
\begin{equation} \label{k-cond}
\lam^k > 2 + \eta^{-1} d_M.
\end{equation}
This will be the $k$ in Theorem~\ref{thm-mld}, and we fix it for the rest of
the proof.

\subsection{Substitution map.}
We are going to define a map
$$
f:\ \Tk \to (\phi^{-k} \Tk)^*
$$
with the following properties:

\smallskip

(S1) 
$$
(T\in \Tk,\ T+g\in \Tk)
\ \Rightarrow\ f(T+g) = f(T)+g;
$$

\smallskip

(S2) $\{f(T):\ T \in \Tk\}$ is a tiling of $\R^d$, that is,
the $\phi^{-k}\Tk$-patches $f(T)$, for $T\in \Tk$,
have supports with disjoint interiors and whose union is all of 
$\R^d$;

\smallskip

(S3) $T,S\in \Tk,\ \phi^{-k} S \in f(T) \ \Rightarrow\ \supp(\phi^{-k} S) \cap
\supp(T) \ne \es$;

\smallskip
                                                                                
(S4) $T,S\in \Tk,\ \supp(\phi^{-k} S) \subset \Int(\supp(T)) \ \Rightarrow\ 
\phi^{-k} S \in f(T)$.

\smallskip 

Observe that the maximal diameter of $\phi^{-k}\Tk$-tiles is not greater than
$\lam^{-k} d_M$, hence (\ref{k-cond}) and (S4) will guarantee that
$f(T)$ is non-empty.

Let $\partial \Tk = \bigcup \{\partial(\supp(T)):\ T\in \Tk\}$.
Note that $\partial\Tk$ is nowhere dense since each tile support is the
closure of its interior.
Let $\{T_1,\ldots,T_m\}$ be a prototile set for $\Tk$. 
Then $\Tk = \bigcup_{i=1}^m (T_i + \Lam_i)$ for some Delone sets $\Lam_i$.
For each $T_i$ choose a ``reference point'' 
\begin{equation} \label{eq-ref}
c(T_i)\in \Int(\supp(T_i)) \setminus \bigcup_{g\in \Lam_i} 
(\phi^k (\partial\Tk) - g),
\end{equation}
which is possible by the Baire Category Theorem. Then define
$c(T_i+g) = c(T_i) + g$ for $g \in \Lam_i$ and let
$c(\Tk) = \{c(T):\ T\in \Tk\}$.
Now for $T\in \Tk$ we define
$$
f(T) = \{\phi^{-k} S:\ \phi^{-k} c(S) \in \Int(\supp(T))\}.
$$
In words, $f(T)$ is the patch of $\phi^{-k}\Tk$-tiles whose reference
points lie in the interior of $\supp(T)$.
Observe that $\phi^{-k} c(\Tk) \cap \partial \Tk = \es$ by (\ref{eq-ref}),
so $f$ is well-defined. Conditions (S2)-(S4) hold by construction.
Since $f(T)$ depends only on the $\phi^{-k}\Tk$ patch of tiles intersecting
$T$, the
condition (S1) holds by (\ref{eq1}). Thus, we have defined the desired
map $f:\ \Tk \to (\phi^{-k} \Tk)^*$.

\subsection{Iterating $\phi^k f$.}
The substitution map $f$, defined above, naturally extends to a map
(also denote by $f$) from
$\Tk^*$ to $(\phi^{-k}\Tk)^*$ as follows: for a patch $P\subset \Tk$ let
$$
f(P) = \bigcup \{f(T):\ T\in P\}.
$$
The properties (S3) and (S4) of the map $f$ imply 
the analog of (S4) for patches:
\begin{equation} \label{S5}
P \subset \Tk,\ S\in \Tk,\ \supp(\phi^{-k}S) \subset \Int(\supp(P))
\ \Rightarrow\ \phi^{-k} S \in f(P).
\end{equation}
Note that $\phi^k f$ maps $\Tk^*$ to itself, so we can iterate it and
obtain the maps $(\phi^k f)^n: \Tk^* \to \Tk^*$.

Now the rough idea for the rest of the proof is that iterating $\phi^k f$ and
rescaling we obtain a self-affine tiling in the limit. More precisely,
one can show that $\supp(\phi^{-kn}(\phi^k f)^n(T))$ converges to a compact
set $A'$ in the Hausdorff metric, for any tile $T \in \Tk$. Then one can prove
that $A'$ is the closure of its interior, and we may consider the tile 
$T' = (A',\ell(T))$. It turns out that $\{T':\ T\in \Tk\}$ is the
desired self-affine tiling. There are a number of technical obstacles
on this route, and we are going to proceed indirectly, utilizing the theory
of substitution Delone sets developed by Lagarias and Wang, with an
extra step coming from \cite{LMS}.

We will need the following simple fact, which is immediate from the
definitions:
for any tiling $\Sk$, compact set $F$, and invertible linear
map $\psi$,
\begin{equation} \label{eq-dumb}
\psi\left([F]^{\Sk} \right) = [\psi F]^{\psi\Sk}.
\end{equation}

\begin{lemma} \label{lem-vspom2}
Let $R\ge \eta$ and $x\in \R^d$. Then
$$
(\phi^k f)([N_R(x)]^{\Tk}) \supset [N_{2R}(\phi^k x)]^{\Tk}.
$$
\end{lemma}

{\em Proof.} Observe that the maximal diameter of a $\phi^{-k}\Tk$-tile
is at most $d_M \lam^{-k}$. Thus, by (\ref{S5}),
$$
f([N_R(x)]^{\Tk}) \supset [N_{R-d_M\lam^{-k}}(x)]^{\phi^{-k} \Tk}.
$$
Then by (\ref{eq-dumb}),
$$
(\phi^k f)([N_R(x)]^{\Tk}) \supset [\phi^k N_{R-d_M\lam^{-k}}(x)]^{\Tk}
\supset [N_{\lam^k R-d_M} (\phi^k x)]^{\Tk}.
$$
Now the desired statement follows since $\lam^k R - d_M > 2R$ for $R\ge \eta$
by (\ref{k-cond}). \qed

\subsection{Substitution Delone set family.} For each prototile $T_j$,
we have a $\phi^{-k}\Tk$ patch $f(T_j)$. Thus, we can write
\begin{equation} \label{sub1}
f(T_j) = \bigcup_{i\le m} (\phi^{-k} T_i + \phi^{-k} D_{ij}),
\end{equation}
where $D_{ij} \subset \Lam_i$ is a finite set. We have by (S2), (S1), and
(\ref{sub1}),
\begin{eqnarray*}
\bigcup_{i\le m}(\phi^{-k}T_i + \phi^{-k} \Lam_i) =
\phi^{-k} \Tk & = & \bigcup_{T\in \Tk} f(T) \nonumber \\
              & = & \bigcup_{j\le m} (f(T_j) + \Lam_j) \nonumber \\
              & = & \bigcup_{j \le m}
\Bigl( \bigcup_{i\le m} (\phi^{-k} T_i + \phi^{-k} D_{ij}) + \Lam_j \Bigr) 
\nonumber \\
 & = & \bigcup_{i\le m} \Bigl(\phi^{-k} T_i + 
\bigcup_{j\le m} (\Lam_j +\phi^{-k}
D_{ij}) \Bigr),
\end{eqnarray*}
where all the unions are essentially disjoint, that is, the tiles and
patches in the unions have disjoint support interiors.
It follows that
\begin{equation} \label{sub2}
\Lam_i = \bigcup_{j \le m}(\phi^k \Lam_j + D_{ij}),\ \ \ i=1,\ldots,m,
\end{equation}
where the unions are disjoint.

Now we need some definitions from \cite{LM,LaWa,LMS}.
A {\em multiset\footnote{Caution: in \cite{LaWa}, 
the word multiset refers to a set with multiplicities.}} 
in $\R^d$ is a
subset $\Xb = X_1 \times \dots \times X_m
\subset (\R^d)^m$
where $X_i \subset \R^d$. We also write
$\Xb = (X_1, \dots, X_m) = (X_i)_{i\le m}$.
Although $\Xb$ is a product of sets, it is convenient to think
of it as a set with types or colors, $i$ being the
color of points in $X_i$. A {\em Delone multiset} is a multiset $\Xb = 
(X_i)_{i\le m}$ where each $X_i$ is  a Delone set.

Consider the following mapping on multisets:
\begin{equation} \label{eq-sub3}
\Phi\left((X_i)_{i\le m}\right) = 
\left(\bigcup_{j\le m} (\phi^k X_j + D_{ij})\right)_{i\le m}
\end{equation}
Then the Delone multiset $\Lb:= (\Lam_i)_{i\le m}$ is a
fixed point of $\Phi$ by (\ref{sub2}), that is,
$\Phi(\Lb) = \Lb$. We say that $\Lb$ is a {\em substitution Delone
multiset} (see \cite[Def.\,3.4]{LMS}). The {\em substitution matrix} is
defined by $\Sf := [|D_{ij}|]_{i,j\le m}$. The substitution Delone set is
said to be {\em primitive} if $\Sf^\ell$ is strictly positive (entry-wise)
for some $\ell\in \Nat$. 

\medskip

\begin{lemma} \label{lem-prim} The substitution Delone multiset
$\Lb:= (\Lam_i)_{i\le m}$ is primitive.
\end{lemma}

{\em Proof.}
By (\ref{sub1}), $\Sf(i,j) = D_{ij}$ counts the number of tiles equivalent to
$T_i$ in $(\phi^k f) (T_j)$. It is easy to see that
$(\Sf)^\ell(i,j)$ counts the number of tiles equivalent to 
$T_i$ in $(\phi^k f)^\ell (T_j)$. Recall that every $\Tk$-tile support 
contains a closed ball of radius $\eta>0$ in its interior.
Thus, $\{T_j\} = [N_\eta(x)]^{\Tk}$ for some $x\in \supp(T_j)$.
Then we can apply Lemma~\ref{lem-vspom2} $\ell$ times to obtain that
\begin{equation} \label{eq-new1}
(\phi^k f)^\ell(T_j) = (\phi^k f)^\ell([N_\eta(x)]^{\Tk}) \supset 
[N_{2^\ell \eta} (\phi^{k\ell} x)]^{\Tk}.
\end{equation}
By the repetitivity of $\Tk$, the right-hand side contains tiles of all
types for $\ell$ sufficiently large, and the claim follows. \qed

\subsection{Self-affine tiling}
For our substitution Delone multiset $\Lb$ there is an {\em adjoint system
of set equations}
\begin{equation} \label{eq-adj}
\phi^k F_j = \bigcup_{i\le m} (F_i + D_{ij}),\ \ \ j\le m.
\end{equation}
It is a well-known fact in the theory of graph-directed iterated function
systems that (\ref{eq-adj}) has a unique solution for which $(F_i)_{i \le m}$
is a family of non-empty compact sets in $\R^d$. It is proved in
Theorems 2.4 and 5.5 of \cite{LaWa} that if $\Lb$ is a primitive
substitution Delone multiset, then each $F_i$ from (\ref{eq-adj}) has
non-empty interior and is the closure of its interior.

\begin{prop} \label{prop-saf}
Let $(F_i)_{i \le m}$ be the solution of  (\ref{eq-adj}). Then
$$
\R^d=\bigcup_{i \le m} (F_i + \Lam_i),
$$
and the sets in the right-hand side have disjoint interiors.
In other words, we obtain a tiling of $\R^d$ by translates of $F_i$.
\end{prop}

In the terminology of \cite{LMS}, the proposition means that the
substitution Delone multiset $(\Lam_i)_{i\le m}$ is {\em representable}
by tiles $(F_i,i)$.

\medskip

{\em Proof.} We are going to use Theorem 3.7 from \cite{LMS}, which in turn
is based on Theorem 7.1 of \cite{LaWa}. 

A {\em $\Lb$-cluster} is a finite multiset $(\Gam_i)_{i\le m}$ where
$\Gam_i\subset \Lam_i$ for all $i$. We have $0 \in \Lam_i$ for all $i$
since $T_i \in \Tk$. For $i \le m$ consider the $\Lb$-cluster
$\be^{(j)} = (e^{(j)}_i)_{i \le m}$ where $e^{(j)}_i = \es$ for $j\ne i$
and $e^{(j)}_j = \{0\}$. 
There is a natural 1-to-1 correspondence between $\Tk$-patches and
$\Lb$-clusters: given a $\Tk$-patch $P$, we can consider the cluster
$\Gamb(P) := (\Gam_i)_{i\le m}$ where 
$P = \bigcup_{i\le m} (T_i + \Gam_i)$, and conversely, every
$\Lb$-cluster arises this way. 

\begin{lemma} \label{lem-new1}
{\bf (i)} For any $\Tk$-patch $P$ we have
$$
\Phi(\Gamb(P)) = \Gamb(\phi^k f(P)).
$$

{\bf (ii)} We have
$$
\Phi^\ell(\be^{(j)}) = \Gamb((\phi^k f)^\ell(T_j)),\ \ \ j\le m,\ \ell \in \Nat.
$$
\end{lemma}

{\em Proof of the lemma.} (i) Let $P = \bigcup_{i\le m} (T_i + \Gam_i)$
and $\Gamb = \Gamb(P) = (\Gam_i)_{i\le m}$. We have, essentially
repeating the calculation which follows (\ref{sub1}),
\begin{eqnarray*}
f(P) & = & \bigcup_{j\le m} (f(T_j) + \Gam_j) \\
& = & \bigcup_{i\le m} \Bigl(\phi^{-k} T_i +
\bigcup_{j\le m} (\Gam_j +\phi^{-k}
D_{ij}) \Bigr).
\end{eqnarray*}
Multiplying by $\phi^k$ and comparing to (\ref{eq-sub3}) yields the
desired equality.

(ii) For $\ell=1$ this holds by (\ref{sub1}) and (\ref{sub2});
for $\ell>1$ this follows from part (i). \qed

\medskip

{\em Conclusion of the proof of Proposition~\ref{prop-saf}.} 
A cluster is said to be {\em legal} if it occurs
in $\Phi^\ell(\be^{(j)})$ for some $j\le m$ and $\ell\in \Nat$.
In Theorem 3.7 of \cite{LMS} it is proved that $\Lb$ is representable by
tiles $(F_i,i)$ if and only if every $\Lb$-cluster is legal.

Let $\Gamb$ be any $\Lb$-cluster. Then $\Gamb = \Gamb(P)$ for some
$\Tk$-patch $P$. Fix any $j \le m$ and $x\in \supp(T_j)$ as in 
(\ref{eq-new1}). Since $\Tk$ is repetitive, (\ref{eq-new1}) implies that
there exists $\ell\in \Nat$ such that $(\phi^k f)^\ell(T_j)$ contains
a $\Tk$-patch equivalent to $P$. By Lemma~\ref{lem-new1}(ii), it
follows that $\Phi^\ell(\be^{(j)})$ contains a $\Lb$-cluster
equivalent to $\Gamb$, and the proof is complete. \qed

\medskip

\begin{prop} \label{prop-new2}
Let $(F_i)_{i \le m}$ be the solution of  (\ref{eq-adj}). For $i\le m$
let $T_i' = (F_i, i)$ and 
$$
\Tk' = \bigcup_{i\le m} (T'_i + \Lam_i).
$$
Then 

{\bf (i)} $\Tk'$ is a self-affine tiling of $\R^d$ with expansion $\phi^k$;

{\bf (ii)} $\Tk'$ is MLD with $\Tk$.
\end{prop}

{\em Proof.} We already proved that $\Tk'$ is a tiling in 
Proposition~\ref{prop-saf}. The fact that it is repetitive can be shown 
directly (as in the argument above dealing with legal clusters), but it
will also follow from part (ii). Then for part (i) it remains to check
the ``geometric substitution'' property, see Definition~\ref{def-ld}.
We define
$$
\om(T'_j + x) = \om(T'_j) + x:= \bigcup_{i\le m} (T'_i + D_{ij}) + \phi^k x,\ 
\ \ \mbox{for}\ x\in \Lam_j.
$$
The right-hand side is a patch with support $\phi^k(\supp(T_j'+x))$ by
(\ref{eq-adj}); and it is a $\Tk'$-patch since $\phi^k x + D_{ij}
\subset \Lam_i$ by (\ref{sub2}). The property (ii) of Definition~\ref{def-ld}
holds by construction.

(ii) Recall that $A_j = \supp(T_j)$ and $F_j = \supp(T_j')$. Let 
$$
C = \max_{j \le m} \rho_H(A_j, F_j),
$$
where $\rho_H$ denotes the Hausdorff metric. 
We claim that $\Tk \ld \Tk'$ with a radius $C$. Indeed,
suppose 
\begin{equation} \label{er1}
[N_C(x)]^{\Tk} = [N_C(y)]^{\Tk} + (x-y),
\end{equation}
and let $T'\in \Tk'$ be such that $\supp(T')$ contains $x$. Then
$T' = T_j' + x$ for some $x\in \Lam_j$ and $T = T_j + x \in \Tk$ satisfies
$$
\rho_H(\supp(T),\supp(T')) \le C \ \Rightarrow\ N_C(\supp(T)) \supset
\supp(T') \ni x.
$$
It follows that $\supp(T) \cap N_C(x) \ne \es$, hence $T \in [N_C(x)]^{\Tk}$.
By (\ref{er1}),
$$
T + (y-x) \in \Tk \ \Rightarrow\ T' + (y-x) \in \Tk'.
$$
This proves that $[x]^{\Tk'} \subset [y]^{\Tk'} + (x-y)$, and the opposite
inclusion is proved reversing the roles of $x$ and $y$. 

In the argument above, we did not use any properties of $\Tk$ and $\Tk'$
other than the fact that there is a 1-to-1 correspondence between their
prototiles and the Delone sets $\Lam_i$ are the same for both tilings. Thus,
it also shows $\Tk' \ld \Tk$ with the same radius $C$. 
This completes the proof of the proposition and the proof of the main theorem,
modulo the proof of Proposition~\ref{prop-reduce} provided in the next
section. \qed


\section{Proof of Proposition~\ref{prop-reduce}}

In the next lemma we collect some elementary properties which will be needed.
We will use the following notation
for a compact set $F \subset \R^d$ and $r > 0$:
$$
F^{-r} := \{x \in F:\, \dist(x,\partial F) \ge r\}.
$$

\begin{lemma}\label{lem-vspom}
{\bf (i)} If $\Tk_1\ld\Tk_2$
with a radius $R$, then for any $L>R$ and any compact set
$F\subset \R^d$,
\[ \ [N_L(F)]^{\Tk_1} = [N_L(F-g)]^{\Tk_1} + g\ \Rightarrow\
[N_{L-R}(F)]^ {\Tk_2} = [N_{L-R}(F-g)]^{\Tk_2} + g.\]
 
{\bf (ii)} If $\Tk_1\ld\Tk_2$
with a radius $R$, then for any compact set
$F\subset \R^d$,
$$
[F]^{\Tk_1} = [F-g]^{\Tk_1} + g \ \Rightarrow\
[F^{-R}]^{\Tk_2} = [(F-g)^{-R}]^{\Tk_2} + g.
$$

{\bf (iii)} If $\Tk_1\ld\Tk_2$
with a radius $R_1$ and $\Tk_2\ld\Tk_3$
with a radius $R_2$, then $\Tk_1 \ld \Tk_3$ with a radius $R_1+R_2$.

{\bf (iv)} If $\psi$ is a linear map such that $|\psi(x)| \ge \gamma|x|$ for
all $x\in \R^d$, for some $\gamma>0$, and $F\subset \R^d$ is a compact set,
then
$$
\psi(F^{-R/\gam}) \subset (\psi F)^{-R}.
$$
\end{lemma}

{\em Proof.} 
(i) is immediate from the definition; the main point is that if 
$\supp(T)$ intersects $N_{L-R}(F)$ for some $T \in \Tk_2$, then 
$T \in [x]^{\Tk_2}$ for $x\in \supp(T) \cap N_{L-R}(F)$, and
$[N_R(x)]^{\Tk_1} \subset [N_L(F)]^{\Tk_1}$.

(ii) follows from (i) since $N_R(F^{-R}) \subset F$.

(iii) follows from (i) as well.

(iv) holds since $\dist(\psi(x), \partial(\psi F)) = \dist(\psi(x),
\psi(\partial F)) \ge \gam \,\dist (x, \partial F)$. \qed

\begin{lemma} \label{lem-ld1}
Suppose that $\Tk$ is a pseudo-self-affine tiling with expansion $\phi$,
such that $\phi \Tk \ld \Tk$ with a radius $R>0$ and $|\phi x| \ge \lam x$
for all $x\in \R^d$. Then for all $\ell \ge 0$ we have that

{\bf (i)}
$\phi^{-\ell} \Tk \ld \phi^{-\ell-1} \Tk$ with a radius $R \lam^{-\ell-1}$;

{\bf (ii)}
$\Tk \ld \phi^{-\ell}\Tk$ with a radius $R(\lam-1)^{-1}$.
\end{lemma}

{\em Proof.} (i) Suppose 
$$
[N_L(x)]^{\phi^{-\ell}\Tk} = [N_L(y)]^{\phi^{-\ell}\Tk} + (x-y)
$$ 
for $L \ge R\lam^{-\ell-1}$.
Then by (\ref{eq-dumb}), 
$$
[\phi^{\ell+1} N_L(x)]^{\phi\Tk} = [\phi^{\ell+1} N_L(y)]^{\phi\Tk} +
\phi^{\ell+1}(x-y),
$$
and Lemma~\ref{lem-vspom}(i) implies
$$
[(\phi^{\ell+1} N_L(x))^{-R}]^{\phi\Tk} = 
[(\phi^{\ell+1} N_L(y))^{-R}]^{\phi\Tk} + \phi^{\ell+1}(x-y).
$$
In view of Lemma~\ref{lem-vspom}(iv), 
$$
\left[\phi^{\ell+1} \left( (N_L(x))^{-R\lam^{-\ell-1}}\right) \right]^{\Tk}
 = \left[\phi^{\ell+1} \left( (N_L(y))^{-R\lam^{-\ell-1}}\right) \right]^{\Tk}
+ \phi^{\ell+1}(x-y) .
$$
Finally, since $(N_L(x))^{-r} = N_{L-r}(x)$, applying (\ref{eq-dumb}) yields
$$
[N_{L - R\lam^{-\ell-1}}(x)]^{\phi^{-\ell-1}\Tk} =
[N_{L - R\lam^{-\ell-1}}(y)]^{\phi^{-\ell-1}\Tk} + (x-y),
$$
as desired. 

(ii) is immediate from part (i) and Lemma~\ref{lem-vspom}(iii). \qed

\medskip

{\em Proof of Proposition~\ref{prop-reduce}.}
We are given a pseudo-self-affine tiling $\Sk$ such that $\phi \Sk \ld \Sk$
with a radius $R$. For every tile $S \in \Sk$, with $A = \supp(S)$ we consider
a new tile
$$
T_S:= \left( A, \langle [N_L(A)]^{\Sk}, S\rangle \right),
$$
where $\langle P, S \rangle$ denotes the equivalence class of a patch $P
\subset \Sk$ with marked tile $S\in P$. The positive number $L$ will be
specified later. The new tiling is
$$
\Tk := \{T_S:\ S \in \Sk\}.
$$
There is a finite number of tile types in $\Tk$ by FLC (though, most likely,
it is huge compared to the number of tile types in $\Sk$).
It is clear that $T_{S'} = T_S + g$ implies $S' = S+g$, so
$\Tk \ld \Sk$ with a radius $0$. By definition, $\Sk \ld \Tk$ with a radius
$L$. It remains to check that for $S,S'\in \Sk$ with supports $A,A'$, 
for any $k\ge 1$,
\begin{equation} \label{equ1}
T_{S'} = T_S + g \ \Rightarrow\ [A']^{\phi^{-k} \Tk} = 
[A]^{\phi^{-k} \Tk} + g.
\end{equation}
Equivalently, we need to verify that for every $T''= T_{S''}\in \Tk$ with
support $A''$ satisfying $\phi^{-k} A'' \cap A' \ne \es$, we have 
$\phi^{-k} T'' + g\in \phi^{-k} \Tk$. What we do know is that
the labels of $T=T_S$ and $T' = T_{S'}$ are the same, hence
$$
[N_L(A')]^{\Sk} = [N_L(A)]^{\Sk} + g.
$$
Thus, by Lemma~\ref{lem-ld1}(ii) and Lemma~\ref{lem-vspom}(i),
\begin{equation} \label{equ1.2}
[N_{L-R(\lam-1)^{-1}}(A')]^{\phi^{-k} \Sk} =
[N_{L-R(\lam-1)^{-1}}(A)]^{\phi^{-k} \Sk} +g.
\end{equation}
Thus, $\phi^{-k} S'' + g \in \phi^{-k} \Sk$, that is,
$\phi^{-k} S'' + g = \phi^{-k} S'''$ for some $S''' \in\Sk$.
Let $A'''$ be the support of $S'''$. We need to show that
$\phi^{-k}T'' + g = \phi^{-k} T'''$, where $T''' = T_{S'''}$. We already
know the equality of supports, so it remains to show that the labels are the
same. This will follow if we verify that
$\phi^{-k}[N_L(A'')]^{\Sk} +  g = \phi^{-k}[N_L(A''')]^{\Sk},$
or equivalently, that
\begin{equation} \label{equ2}
[\phi^{-k} N_L(A''')]^{\phi^{-k} \Sk}  =  
[\phi^{-k} N_L(A'')]^{\phi^{-k} \Sk} + g.
\end{equation}
However,
$$
[\phi^{-k} N_L(A'')]^{\phi^{-k} \Sk}\subset 
[N_{L\lam^{- k}} (\phi^{-k} A'')]^{\phi^{-k} \Sk}\subset
[N_{(L+d_M)\lam^{-k}} ( A')]^{\phi^{-k} \Sk},
$$
since $d_M(\phi^{-k}\Sk) \le \lam^{-k} d_M$,
and (\ref{equ2}) follows from (\ref{equ1.2}), provided
$L-R(\lam-1)^{-1} > (L+d_M)\lam^{-1} > (L+d_M)\lam^{-k}$. 
Thus, it is enough to choose
$$
L > R\lam(\lam-1)^{-2} + d_M (\lam-1)^{-1},
$$
and the proof is complete. \qed

\medskip


\noindent {\bf Acknowledgment.} I am grateful to Jeong-Yup Lee and
Lorenzo Sadun for helpful discussions.
 
 
\bibliographystyle{amsplain}

\begin{thebibliography}{99}
                                                                                
\bibitem{BS} M. Baake and  M. Schlottmann, Geometric aspects
of tilings and equivalence concepts, in {\em Proc.\ ICQ5}, World Scientific,
Singapore, 1995, pp. 15--21.
                                                                                
\bibitem{Baake-ld} M. Baake, M. Schlottmann, and P. D. Jarvis,
Quasiperiodic tilings with tenfold symmetry and equivalence with
respect to local derivability, {\em J.\ Phys.\ A} {\bf 24} (1991), 4637--54.

\bibitem{HRS} C. Holton, C. Radin, and L. Sadun,
Conjugacies for tiling dynamical systems, {\em Preprint} math.DS/0307259.

\bibitem{KV} R. Kenyon and A. Vershik,
Arithmetic construction of sofic partitions of hyperbolic toral automorphisms, 
{\em  Ergodic Theory Dynam.\ Systems} {\bf 18}  (1998),  no.\ 2, 357--372.

\bibitem{LaWa} J. C. Lagarias and Y. Wang,
Substitution Delone sets, {\em  Discrete Comput.\ Geom.} {\bf 29} (2003),
175--209.

\bibitem{LMS} J.-Y. Lee, R. V. Moody, and B. Solomyak,
Consequences of pure point diffraction spectra for multiset 
substitution systems, {\em  Discrete Comput.\ Geom.} {\bf  29}  (2003), 
525--560.

\bibitem{LM} D. Lind and B. Marcus, {\em An introduction to symbolic
dynamics and coding}, Cambridge University Press, Cambridge, 1995.

\bibitem{petersen} K. Petersen, Factor maps between tiling dynamical systems,
{\em Forum Math.} {\bf 11} (1999), 503--512.
                                                                                
\bibitem{Thesis} N. M. Priebe, Detecting hierarchy in tiling dynamical
systems via derived Vorono\"{\i} tesselations,
Ph.D. Thesis, University of North Carolina
at Chapel Hill, 1997.

\bibitem{PF} N. M. Priebe,
Towards a characterization of self-similar tilings in terms of derived
Vorono\"{\i} tesselations, {\em Geom.\ Dedicata} {\bf 79} (2000), 239--265.

\bibitem{PrSo} N. M. Priebe and B. Solomyak,
Characterization of planar pseudo-self-similar tilings,
{\em Discrete Comput.\ Geom.} {\bf 26} (2001), 289--236.

\bibitem{rasad} C. Radin and L. Sadun, Isomorphism of hierarchical structures,
{\em Ergodic Theory Dynam.\ Systems},  {\bf 21}  (2001), 1239--1248.

\bibitem{rand} B. Rand, Equivalence of self-similar and pseudo-self-similar
tiling spaces in $\R^2$, Preprint, 2004.

\bibitem{Robi} E. A. Robinson, Jr., 
Symbolic dynamics and tilings of $\R^d$,  in {\em 
Symbolic dynamics and its applications},  
{\em Proc.\ Sympos.\ Appl.\ Math.}, {\bf 60}, 
Amer.\ Math.\ Soc., Providence, RI, 2004, pp. 81--119.

\bibitem{Rud} D. J. Rudolph, 
Markov tilings of $\R^n$ and representations of $\R^n$ actions,
in {\em Measure and measurable dynamics (Rochester, NY, 1987)}, 
Amer.\ Math.\ Soc., Providence, RI, 1989, pp. 271--290.

\bibitem{SiVer} N. Sidorov and A. Vershik,
Bijective arithmetic codings of hyperbolic automorphisms of the $2$-torus, 
and binary quadratic forms,  
{\em J.\ Dynam.\ Control Systems} {\bf 4}  (1998),  no.\ 3, 365--399.

\bibitem{SolTil} B. Solomyak,
Dynamics of Self-Similar Tilings,
{\em Ergodic Theory Dynam.\ Systems} {\bf 17} (1997), 695--738.
Corrections, {\em Ergodic Theory Dynam.\ Systems} {\bf 19} (1999), 1685.

\bibitem{Ver} A. Vershik,
 Arithmetic isomorphism of hyperbolic automorphisms of a torus and of sofic 
shifts, {\em  Funct.\ Anal.\ Appl.} {\bf 26}  (1992),  no.\ 3, 170--173.

\end{thebibliography}

\end{document}